\documentclass[letterpaper, 10pt, conference]{ieeeconf}

\IEEEoverridecommandlockouts   
\usepackage{amsmath,amssymb,amsthm,amsfonts}
\usepackage{dsfont}
\usepackage{booktabs}
\usepackage{graphicx}
\usepackage{xcolor}
\usepackage{mathtools}
\usepackage{array}

\usepackage[noadjust]{cite}
\usepackage{hyperref}
\hypersetup{colorlinks=true,linkcolor=blue,citecolor=blue,urlcolor=cyan}

\newcommand{\E}{\mathbb{E}}
\newcommand{\R}{\mathbb{R}}
\newcommand{\Prob}{\mathbb{P}}
\newcommand{\tr}{\operatorname{tr}}

\DeclareMathOperator*{\argmin}{arg\,min}
\newcommand{\dd}{\mathrm{d}}
\newcommand{\pos}[1]{(#1)^{+}}

\newtheorem{theorem}{Theorem}

\theoremstyle{definition}

\newtheorem{remark}[theorem]{Remark}

\title{\LARGE \bf Constrained Deep Inventory Management Using Forward-Backward SDEs}

\author{Keunwoo Lim and Michael Bloem%
\thanks{Keunwoo Lim is with the Department of Statistics, University of Washington,
Seattle, WA, USA. Michael Bloem is with Amazon, Bellevue, WA, USA.
Emails: kwlim@uw.edu, bloemm@amazon.com.}}

\begin{document}
\maketitle
\thispagestyle{empty}
\pagestyle{empty}

\begin{abstract}
We study multi-product replenishment subject to aggregate operating limits imposed pathwise in place of expectation. On-hand and
in-transit inventory, warehouse capacity and an emission allowance define a safe set, and the resulting multi-product lost-sales problem is posed as a safe stochastic optimal control problem. The Hamilton--Jacobi--Bellman (HJB) equation is represented as a system of forward--backward stochastic differential equations (FBSDEs) and solved by learning the value function
and its gradient, with a quadratic program derived by control barrier function (CBF) is solved at each decision epoch. Because purchase orders enter the in-transit pipeline only, the in-transit constraints
have relative degree one while the on-hand constraints have relative degree two. We extend the stochastic high-relative-degree construction of \cite{clark2021control} from a single barrier to four coupled aggregate constraints of mixed relative degree, and we establish its initial-condition requirement by design. 
\end{abstract}

\section{Introduction}

Inventory replenishment across a product assortment is able to be formulated as a stochastic control problem in which decisions for individual products interact through shared operating limits. A retailer reviews each product periodically, observes its inventory position, and places a purchase
order under uncertain demand, an expensive stockout, and a holding cost on stock that arrives too early. At assortment scale the problem does not decompose, because the products draw on a common warehouse of finite size and increasingly a shared carbon allowance. Each of these limits applies to a weighted sum across the assortment, so decisions that are individually optimal may be jointly infeasible.

Two features distinguish this setting from its single-product counterpart. First, the state dimension grows with the number of products, so any mechanism enforcing the shared limits must remain tractable as the assortment grows. Second, a purchase order does not act on the constrained quantity directly: ordered units enter an in-transit pipeline and reach the shelf after a vendor lead time.
In the notion of control theory, these constraints have different relative degree.

Several lines of work have reference to this problem. Capacitated inventory theory establishes that a modified base-stock policy for period-dependent targets is optimal for a relaxed multi-item model with random and seasonally varying demand \cite{aviv2001capacitated}, and its continuous-time counterpart yields an $(s,S)$ policy for a Brownian production and inventory system with finite production capacity \cite{wu14optimal}; in both the limit
constrains the quantity ordered or produced, which is a constraint on the control rather than on the state. At assortment scale, policies learned by differentiating a simulator driven by historical data represent the current practice, and aggregate capacity has been addressed there by learning a neural coordinator \cite{eisenach2024neural}; the limits are then respected in expectation, with residual violations along realized trajectories measured rather than eliminated. Section~\ref{sec:related}
reviews these literature in full.

We take the alternative route of enforcing the shared limits pathwise. Control barrier functions \cite{ames2016control} certify forward invariance of a
safe set through a condition imposed pointwise on the control, and
\cite{clark2021control} develops the stochastic counterpart. Combining the probabilistic representation of the HJB equation \cite{exarchos2018stochastic} with function approximation yields the deep FBSDE controller
\cite{pereira2019learning}, and \cite{pereira2021safe} enforce safety within it through a differentiable convex optimization layer \cite{amos2017optnet}. In contrast, the on-hand limits considered here do not have relative degree one, and removing that restriction is the technical contribution of Section~\ref{sec:cbf}: we extend the high-relative-degree construction from a single barrier to four coupled aggregate constraints of mixed relative degree, and solve them together with $2M$ per-product bounds as
one quadratic program per sampled trajectory, where $M$ is the number of products.

Operating limits in inventory control are rarely static. Warehouse space and receiving capacity are contracted and renegotiated periodically, and carbon allowances are issued per compliance period. When the parameters of a capacitated inventory problem vary with the period, the optimal target becomes period-dependent \cite{kapuscinski1998capacitated,aviv2001capacitated}, so a policy specified by a target level must be recomputed for each phase. 
A piecewise-constant limit is also the case a barrier construction accommodates without modification, since the safe set is fixed between changes and no
derivative of the limit enters the conditions.

The contribution of this paper is as follows. We formulate multi-product replenishment with aggregate state constraints on on-hand and in-transit inventory as a safe stochastic optimal control problem in which the limits hold pathwise (Section~\ref{sec:model}). We extend the stochastic high-relative-degree barrier construction of \cite[Section 4.4]{clark2021control} from a single barrier to four coupled constraints of mixed relative degree, and give sufficient conditions under which the required initial condition holds by design (Section~\ref{sec:cbf}). We then show that the guarantee extends to budgets following a piecewise-constant schedule (Theorem~\ref{thm:pwc}), and report a scaling analysis showing that its performance is retained (Section~\ref{sec:exp}).

\subsection{Related work}\label{sec:related}

Shared limits have been studied in inventory theory since
 \cite{evans1967inventory}, with one line recovering tractability by relaxing the shared constraint and coordinating products through a price on the shared resource \cite{rosenblatt1990single,maloney1993constrained}. In the stochastic dynamic setting, \cite{federgruen1986inventoryI,federgruen1986inventoryII} established that a modified base-stock policy works under capacity limit for a single product, and \cite{aviv2001capacitated} extended this to capacitated multi-item systems with random and seasonally fluctuating demand, where modified base-stock for period-dependent targets is optimal for a relaxed
model. The continuous-time counterpart is due to \cite{wu14optimal}, who model cumulative production and demand as a two-dimensional Brownian motion with finite production capacity and show that an $(s,S)$ policy is optimal.
In most of this work, the capacity bounds the quantity ordered or produced within a period, whereas the limits considered in our work bound inventory states.

Emission allowances are the second family of aggregate limit and have been treated separately. The paper \cite{benjaafar2012carbon} show how a carbon
footprint enters standard inventory models either as a constraint or as a cost term, and \cite{absi2013lot} distinguish periodic, cumulative, global and rolling formulations of the same budget, whose placement affects the solution; \cite{helmrich2015economic} analyze the resulting complexity. This literature places the limit on quantities produced or transported and is largely deterministic \cite{chen2010environmental,sebatjane2025sustainable}.

Enforcing a limit along every trajectory instead brings in the safe control literature. Control barrier functions \cite{ames2016control} certify forward
invariance of a safe set and reduce enforcement to a pointwise convex constraint. The paper \cite{clark2021control} develops the stochastic theory, including a recursion for barriers of relative degree greater than one, close to \cite{sarkar2020high}; \cite{xiao2021high} give the deterministic high-order construction. Safety has been combined with deep forward--backward representations of the HJB equation by \cite{pereira2021safe}, using a differentiable optimization layer \cite{amos2017optnet,diamond2016cvxpy}.

\textbf{Notation:} Throughout, $v^i$ denotes the $i$-th element of vector $v$, $\nabla$ denotes the gradient with respect to the state $Z$, $\mathrm{diag}(v)$
the diagonal matrix with real vector $v$ on its diagonal, $e_j\in\R^M$ the $j$th standard basis vector, and $\pos{x}=\max\{x,0\}$ with a real vector $x$  applied componentwise. Inequalities between vectors are
componentwise.

\section{Constrained inventory management framework}\label{sec:model}

This section sets up the model: the continuous-time dynamics of on-hand and in-transit inventory for $M$ products, and the cost functional to be minimized.

\subsection{Inventory dynamics}

We consider procurement of $M$ products over a finite horizon $[0,T]$. Write $X_t^i$ and $L_t^i$ for the on-hand and in-transit inventory of product $i$, collected into $X_t,L_t\in\R^M$ and $Z_t=(X_t^\top,L_t^\top)^\top\in\R^{2M}$. The control $u_t\in\R^M$ is a vector of order rates, so $\int_{t_1}^{t_2}u_s^i\,\dd s$ is the quantity of product $i$ ordered over $[t_1,t_2]$.  The exogenous process $Q_t=(d_t^\top, (c^{p}_t)^\top, (c^{o}_t)^\top, (c^{h}_t)^\top, g_t^\top)^\top\in\R^{5M}$ collects the demand rate $d_t$, unit stockout penalty $c^{p}_t$, unit order cost
$c^{o}_t$, unit holding cost $c^{h}_t$, and inbound throughput limit $g_t$. It evolves independently of states and actions, which is the exogenous-input structure of \cite{sinclair2023hindsight,maggiar2025structure}. The state evolves according to
\begin{align}\label{eq:dynamics}
  &\dd Z_t = \bigl(\,f(Z_t,Q_t) + G u_t\,\bigr)\dd t + \Sigma\,\dd W_t,\\
  &f = \begin{bmatrix}\, b_{\mathrm{in}}(L_t,g_t) - b_{\mathrm{out}}(X_t,d_t)\, \\ \,-\,b_{\mathrm{in}}(L_t,g_t)\,\end{bmatrix},
  \,\,
  G = \begin{bmatrix} 0 \\ I_M \end{bmatrix},
  \notag
\end{align}
where $W_t$ is a $2M$-dimensional Brownian motion and diffusion
$\Sigma=\mathrm{diag}(\sigma_X I_M,\ \sigma_L I_M)$ with scalars $\sigma_X,\sigma_L>0$ and identity matrix $I_M$. The two drift functions are
\begin{align}\label{eq:flows}
  b_{\mathrm{in}}(L,g) = \min\Bigl(\,\frac{L^{+}}{\tau_{\mathrm{in}}},\, g\,\Bigr),
  b_{\mathrm{out}}(X,d) = \min\Bigl(\,\frac{X^{+}}{\tau_{\mathrm{out}}},\, d\,\Bigr),
\end{align}
applied componentwise, with $(\cdot)^{+}$ the positive part.

\begin{remark}[Interpretation of $\tau_{\mathrm{in}}$ and $\tau_{\mathrm{out}}$]
\label{rem:tau}
The constant $\tau_{\mathrm{in}}$ represents the mean transit time of the pipeline: since $b_{\mathrm{in}}=L/\tau_{\mathrm{in}}$ in the uncongested regime, a delay of unit in transit follows the exponential distribution with scale $\tau_{\mathrm{in}}$. This differs
from the fixed integer lead time of the discrete lost-sales literature
\cite{zipkin2008old}, under which every unit ordered at $t$ arrives exactly at $t+\ell$, where $\ell$ is the deterministic lead time. Diffusion based inventory dynamics also have a history in inventory control, where
cumulative demand and supply are modeled directly as Brownian motions \cite{wu14optimal}.
Analogously, $\tau_{\mathrm{out}}$ is the time constant of the outflow, and
$\min(X/\tau_{\mathrm{out}},d)$ relaxes the reflecting lost-sales boundary.
\end{remark}

\subsection{Cost}

For a Markov control $\upsilon=(u_t)$ in the admissible class $\mathcal{U}$ of uniformly bounded, Lipschitz feedback laws, the cost functional is
\begin{align}\label{eq:cost}
\begin{split}
  J^{\upsilon}(Z,t) = \E\,\Bigl[\,&\int_t^T \Bigl(\, q(Z_s,Q_s) + (c^{o}_s)^\top u_s\\
    &+ \tfrac12 u_s^\top R u_s \,\Bigr)\,\dd s + \phi(Z_T) \,\Big|\, Z_t = Z\,\Bigr],
\end{split}
\end{align}
with quadratic regularization parameter $R\in \mathbb{R}^{M \times M}$ and
\begin{align*}
  q(Z,Q) \;=\; (c^{h})^\top X^{+} \;+\; (c^{p})^\top\Bigl(\,d - \frac{X^{+}}{\tau_{\mathrm{out}}}\,\Bigr)^{\!+}.
\end{align*}
The first term is holding cost and the second term is the underage cost, aligned with the outflow in \eqref{eq:flows}. The terminal cost values leftover stock at its purchase cost net of salvage at a fraction
$\varsigma\in[0,1)$ of price,
\begin{align*}
  \phi(Z_T) = {\pos{c^{o}_T - \varsigma\, c^{p}_T}}^\top (X_T + L_T),
\end{align*}
so that leftover inventory is a net loss and the controller has reason to wind down as
the horizon approaches. The corresponding value function is $V(Z,t)=\inf_{\upsilon\in\mathcal{U}}J^{\upsilon}(Z,t)$, satisfying
the HJB equation
\begin{align}\label{eq:hjb}
\begin{split}
  &\frac{\partial V}{\partial t}
  + \inf_{u}\Bigl[\tfrac12\tr\Bigl(\frac{\partial^2 V}{\partial Z^2}\Sigma\Sigma^\top\Bigr)
  + \Bigl(\frac{\partial V}{\partial Z}\Bigr)^{\!\top}\!\bigl(f + Gu\bigr)\\
  &\qquad
  + q + (c^{o}_t)^\top u + \tfrac12 u^\top R u\Bigr] = 0,
\end{split}
\end{align}
with $V(Z,T)=\phi(Z)$.

\section{Safe set and control barrier functions}\label{sec:cbf}

This section constructs the pointwise conditions that keep the state in the safe set. We first define the safe set and compute the relative degree of each barrier, which separates them into two families. For the family the control reaches directly we recover a first-order condition (Theorem~\ref{thm:safety1}); for the family it reaches only through the pipeline we establish a coupled second-order condition (Theorem~\ref{thm:rd2}) and give a sufficient condition.
\subsection{The safe set}

Three families of constraint are natural here: validity constraints, aggregate capacity
\cite[Section 4.1]{eisenach2024neural}, and an emission budget
\cite[Section 5]{chen2010environmental}. Inventories and orders are nonnegative,
\begin{align}\label{eq:basic}
  X_t \ge 0, \qquad L_t \ge 0, \qquad 0 \le u_t \le \tilde{u},
\end{align}
where $\tilde{u}$ is a wide per-product bound on the order rate discussed in
Section~\ref{sec:layer}. Capacity is measured in shelf space, and we set every unit to occupy the same space, giving uniform capacity weights. Emissions per unit vary across products, since a heavy or air-freighted item emits considerably more than a small light one, which leads that the emission weights are heterogeneous. With capacity weights $w_X,w_L$ and emission
weights $w_{E_X},w_{E_L}$, the four aggregate budgets are
\begin{align}\label{eq:agg}
\begin{split}
  &w_X^\top X_t \le H^X, w_L^\top L_t \le H^L,\\ 
  &w_{E_X}^\top X_t \le H^{E_X}, w_{E_L}^\top L_t \le H^{E_L}.
\end{split}
\end{align}
The safe set is $\mathcal{C}=\bigcap_{j=1}^{2M+4}\mathcal{C}_j$ with
$\mathcal{C}_j=\{Z: h_j(Z)\ge 0\}$ and
\begin{align}\label{eq:total_barriers}
\begin{split}
  &h_j(Z) = X^j, \, h_{j+M}(Z) = L^j, \, j = 1,\dots,M, \\
  &h_{2M+1} = H^X - w_X^\top X, \,
  h_{2M+2} = H^L - w_L^\top L, \\
  &h_{2M+3} = H^{E_X} - w_{E_X}^\top X, \,
  h_{2M+4} = H^{E_L} - w_{E_L}^\top L.
\end{split}
\end{align}

Two properties of \eqref{eq:agg} should be stated plainly because they shape what can be claimed. The budgets are aggregate: with uniform weights, $w_X^\top X\le H^X$ cannot distinguish $H^X$ units spread evenly across the assortment from $M H^X$ units concentrated on one product. In addition, they constrain the state, not the control. The classical capacitated results, including \cite{federgruen1986inventoryI, aviv2001capacitated}, bound the order or production quantity, which is a constraint on $u$; \eqref{eq:agg} bounds stored inventory, which the control reaches only through the pipeline.

Because $G=[0;I_M]^\top$ has a zero on-hand block, the control enters the barriers very differently. For an in-transit barrier, differentiating once already exposes $u$: $\partial h_{j+M}/\partial Z = e_{j+M}$ gives
$(\partial h_{j+M}/\partial Z)^\top G u = u^j$. For an on-hand barrier,
$(\partial h_j/\partial Z)^\top G = 0$ holds, so the first-order condition is vacuous and one more differentiation is required. The in-transit constraints therefore have relative degree one and the on-hand constraints relative degree two.

\subsection{In-transit conditions (relative degree one)}\label{sec:cbf1}

We first record the standard stochastic CBF statement, which we use for the in-transit
rows. It is based on \cite[Theorem 1]{pereira2021safe} restated for our barriers and exogenous
process.

\begin{theorem}[Safety, relative degree one]\label{thm:safety1}
Let $\alpha:\R\to\R$ be continuous and increasing with $\alpha(0)=0$, and let $h$ be twice continously differentiable with $Z_0$ satisfying $h(Z_0)\ge 0$. If, along the dynamics \eqref{eq:dynamics}
under a control $(u_t)\in\mathcal{U}$,
\begin{align*}
\begin{split}
  \Bigl(\frac{\partial h}{\partial Z}(Z_t)\Bigr)^{\!\top}&\!\bigl(f(Z_t,Q_t)+Gu_t\bigr)+ \tfrac12\tr\Bigl(\frac{\partial^2 h}{\partial Z^2}(Z_t)\,\Sigma\Sigma^\top\Bigr)\\
  & \ge -\alpha\bigl(h(Z_t)\bigr),
\end{split}
\end{align*}
for all $t\in[0,T]$, then $h(Z_t)\ge 0$ for all $t\in[0,T]$ with probability one.
\end{theorem}

For $\gamma>0$, in the case of $\alpha(h)=\gamma h$ to $h_{j+M}=L^j$ gives, since the
barrier is linear so its Hessian vanishes,
\begin{align}\label{eq:lbL}
  u^j \;\ge\; b_{\mathrm{in}}^j(L,g) - \gamma L^j \;=:\; \underline u_L^j .
\end{align}
For the aggregate in-transit rows, we obtain
\begin{align}
  &w_L^\top u \;\le\; w_L^\top b_{\mathrm{in}}(L,g) + \gamma\bigl(H^L - w_L^\top L\bigr),\label{eq:rowL}\\
  &w_{E_L}^\top u \;\le\; w_{E_L}^\top b_{\mathrm{in}}(L,g) + \gamma\bigl(H^{E_L} - w_{E_L}^\top L\bigr).\label{eq:rowL_2}
\end{align}
\subsection{Second-order conditions (relative degree two)}\label{sec:rd2}

For the on-hand rows we use the high-relative-degree construction of
\cite[Section 4.4]{clark2021control}, close to the stochastic construction of \cite{sarkar2020high}. We denote this result from \cite[Theorem 4]{clark2021control} in Theorem \ref{thm:rd2}. Note that Theorem \ref{thm:rd2} requires the twice continuous differentiability of $f$, which does not hold for the drift in \eqref{eq:dynamics}. Therefore, we assume that $f$ is sufficiently approximated to satisfy the twice continuous differentiability in the theorem, and introduce the smoothing method of $f$,   $b_{\mathrm{in}}$, and $b_{\mathrm{out}}$ in Section \ref{sec:smoothing}.

\begin{theorem}[Safety under coupled constraints of mixed relative degree]\label{thm:rd2}
Let $\Sigma$ be constant, and assume that $f$ is sufficiently approximated to be twice continuously differentiable. Partition the barriers of \eqref{eq:total_barriers} into
\begin{align*}
  \mathcal{J}_1 &= \{\,h_j:\ \nabla h_j^\top G \neq 0,\, 1\leq j \leq 2M+4\,\}, \\
  \mathcal{J}_2 &= \{\,h_j:\ \nabla h_j^\top G = 0, \, 1\leq j \leq 2M+4\,\},
\end{align*}
with each $h_j$ affine. For $h_j\in\mathcal{J}_2$, define $h_{(1, j)} = \nabla h_j^\top f + \gamma h_j$.
Suppose the initial state satisfies
\begin{align}\label{eq:initcond}
\begin{split}
  h_j(Z_0)\ge 0 \text{ for all }h_j\text{ and }
  h_{(1, j)}(Z_0)\ge0
  \text{ for all } h_j\in\mathcal{J}_2,
\end{split}
\end{align}
and that the control $(u_t)\in\mathcal{U}$ satisfies, for all $t\in[0,T]$,
\begin{align}
  \nabla h_j^\top\!\bigl(f+Gu_t\bigr) \;\ge\; -\gamma\, h_j, \label{eq:cond1}
\end{align}  
for every $h_j\in\mathcal{J}_1$, and
\begin{align}
  \nabla h_{(1, j)}^\top\!\bigl(f+Gu_t\bigr)  + \tfrac12\tr\bigl(\nabla^2 h_{(1, j)}\Sigma\Sigma^\top\bigr)
    \ge -\gamma h_{(1, j)}, \label{eq:cond2}
\end{align}
for every $h_j\in\mathcal{J}_2$. Then $\Prob\bigl(\,Z_t\in\mathcal{C},\  \forall t\in[0,T]\,\bigr)=1$.
\end{theorem}

For smoothed, twice continuously differentiable $b_{\mathrm{in}}$ and $b_{\mathrm{out}}$, write $\beta_j = \partial b_{\mathrm{in}}^j/\partial L^j$ and
$a_j = \partial b_{\mathrm{out}}^j/\partial X^j$ for the first derivatives of the two flows,
$\beta'_j = \partial \beta_j/\partial L^j$ and $a'_j = \partial a_j/\partial X^j$ for the corresponding second derivatives, and
$f_X = b_{\mathrm{in}} - b_{\mathrm{out}}$, $f_L = -b_{\mathrm{in}}$ for the two blocks of
the smoothed drift. First derivatives $\beta_j$ represents the control
authority, and $a_j$ represents the marginal service rate. Then, with $h_{(0, j)}=X^j$ linear, the corresponding CBF \eqref{eq:cond2} is
\begin{align}\label{eq:lbX}
\begin{split}
  u^j \ge &-\gamma\bigl(f_X^j+\gamma X^j\bigr)/\beta_j
    - \bigl[(\gamma-a_j)f_X^j + \beta_j f_L^j\bigr]/\beta_j\\
    &- \bigl(\sigma_L^2 \beta'_j - \sigma_X^2 a'_j\bigr) / 2\beta_j
  \coloneqq \underline u_X^j,
\end{split}
\end{align}
by direct computation. Also, for $h_{(0, 2M+1)}=H^X-w_{X}^\top X$, we get $h_{(1, 2M+1)} = -w_X^\top f_X + \gamma\,(H^X-w_X^\top X)$ and
\begin{align}\label{eq:rowX}
\begin{split}
  &\sum_j w_X^j\beta_j\, u^j \le\;
  \sum_j w_X^j(a_j-\gamma)f_X^j- \sum_j w_X^j\beta_j f_L^j\\
  &\,\,+ \tfrac12\Bigl(\sigma_X^2\sum_j w_X^j a'_j - \sigma_L^2\sum_j w_X^j \beta'_j\Bigr)+ \gamma h_{(1, 2M+1)},
\end{split}
\end{align}
and similarly with $h_{(1, 2M+3)} = -w_{E_X}^\top f_X + \gamma\,(H^{E_X}-w_{E_X}^\top X)$,
\begin{align}\label{eq:rowX_2}
\begin{split}
  &\sum_j  w_{E_X}^j\beta_j\, u^j \le
  \sum_j w_{E_X}^j(a_j-\gamma)f_X^j- \sum_j w_{E_X}^j\beta_j f_L^j\\
  &\,\,+ \tfrac12\Bigl(\sigma_X^2\sum_j w_{E_X}^j a'_j - \sigma_L^2\sum_j w_{E_X}^j \beta'_j\Bigr)+ \gamma h_{(1, 2M+3)}.
\end{split}
\end{align}

\begin{remark}[Feasibility of the initial condition \eqref{eq:initcond}]\label{rem:C1}
Condition \eqref{eq:initcond} is not implied by
$Z_0\in\mathcal{C}$ and must be established. Taking $\gamma=1/\tau_{\mathrm{out}}$, which we do
throughout the implementation, discharges it for the per-product on-hand barriers under certain conditions. For the
aggregate on-hand rows,
\begin{align}\label{eq:aggC1}
\begin{split}
  h_{(1, 2M+1)} = {}&w_X^\top b_{\mathrm{out}} - \gamma\,w_X^\top X+ \gamma H^X - w_X^\top b_{\mathrm{in}} ,
\end{split}
\end{align}
and when every product lies below its service threshold, this reduces to
$\gamma H^X - w_X^\top b_{\mathrm{in}}$. Since $w_X = w_L$ and $w_L^\top b_{\mathrm{in}} \le H_L/\tau_{\mathrm{in}}$ on the safe set, a sufficient condition for satisfying \eqref{eq:aggC1} is
\begin{align}\label{eq:sizing}
  H^L\,\tau_{\mathrm{out}} \;\le\; H^X\,\tau_{\mathrm{in}} .
\end{align}
Condition \eqref{eq:sizing} admits an operational reading: $H^L/\tau_{\mathrm{in}}$ is the largest
sustained arrival rate when the pipeline sits at its budget and $H^X/\tau_{\mathrm{out}}$ is the
largest sustained shipping rate when the shelf sits at its budget, so it represents that the inbound
operation cannot deliver faster than the warehouse can clear.
\end{remark}

\subsection{Piecewise-constant budgets}\label{sec:pwc}

The budgets treated in previous sections are constant. In practice they are not: warehouse and receiving
capacity are contracted and renegotiated at known dates, carbon allowances are issued per
compliance period \cite{absi2013lot}, and capacity is itself managed over a horizon
\cite{angelus2002simultaneous}. When the parameters of a capacitated inventory problem vary by
period, the optimal target becomes period-dependent
\cite{kapuscinski1998capacitated,aviv2001capacitated}, and learned formulations that sample
capacity paths in a Haar basis face the same structure \cite{eisenach2024neural}. We therefore
take $H^X, H^L, H^{E_X}, H^{E_L}$ to be piecewise constant, which is the case in which the construction of
Sections~\ref{sec:cbf1} and \ref{sec:rd2} extends with no modification.

Let $0=t_0<t_1<\dots<t_K=T$ and let $H^X, H^L, H^{E_X}, H^{E_L}$ be right-continuous function of time, with $H^X\equiv H^{X(k)}, H^L\equiv H^{L(k)}, H^{E_X}\equiv H^{E_X(k)}, H^{E_L}\equiv H^{E_L(k)}$ on
$[t_k,t_{k+1})$. Write $\mathcal{C}^{(k)}$ and $h_j^{(k)}$ as the safe set and CBF formed with $H^{X(k)}, H^{L(k)}, H^{E_X(k)}, H^{E_L(k)}$, and
\begin{align*}
\begin{split}
  \bar{\mathcal{C}}^{(k)} = \mathcal{C}^{(k)} \cap
  \bigl\{Z: h^{(k)}_{(1, j)}(Z)\ge0 \text{ for every } h_j\in\mathcal{J}_2\bigr\},
\end{split}
\end{align*}
for the states meeting both requirements of \eqref{eq:initcond} at that budget. The barriers
$X^j$ and $L^j$ do not involve $H$, so a change of budget affects only the four aggregate
barriers, and the argument is confined to them. Next theorem supports our final resulting model.

\begin{theorem}[Safety under a piecewise-constant schedule]\label{thm:pwc}
Let $H$ be right-continuous and piecewise constant on the partition
$0=t_0<\dots<t_K=T$, and suppose that on each interval $[t_k,t_{k+1})$ the control satisfies
\eqref{eq:cond1} and \eqref{eq:cond2} with the budget $H^{(k)}$ in force. Assuming $Z_k \in \bar{\mathcal{C}}^{(k)}$ for all $0\leq k \leq K$, we have
$\Prob\bigl(Z_t\in\mathcal{C}\bigl(H(t)\bigr)\ \ \forall t\in[0,T]\bigr)=1$.
\end{theorem}

\subsection{FBSDE formulation}

The constrained Hamiltonian minimization defines the control
\begin{align}\label{eq:qp}
\begin{split}
  &u_t^{\ast}(Z) \in \argmin_{u}\ \Bigl(\frac{\partial V}{\partial Z}\Bigr)^{\!\top} G u+ (c^{o}_t)^\top u + \tfrac12 u^\top R u\\
  &\qquad \text{s.t. \eqref{eq:basic}, \eqref{eq:lbL}, \eqref{eq:rowL}, \eqref{eq:rowL_2}, \eqref{eq:lbX}, \eqref{eq:rowX}, \eqref{eq:rowX_2}},
\end{split}
\end{align}
and the modified HJB equation is \eqref{eq:hjb} with the infimum replaced by evaluation at
$u_t^{\ast}$. Applying the nonlinear Feynman-Kac lemma following
\cite{exarchos2018stochastic} gives the FBSDE system used for computation in Theorem \ref{thm:fbsde}.

\begin{theorem}\label{thm:fbsde}
Assume $u_t^{\ast}(Z)$ and $K_t(Z):\R^{2M}\to\R^{2M}$ are uniformly bounded and Lipschitz
in $Z$, and that the solution $V$ of the modified HJB equation exists with
$\partial V/\partial t$ and $\partial^2V/\partial Z^2$ continuous. Then $V$ satisfies
\begin{align*}
  \begin{cases}
    \dd Z_t = \bigl(f + Gu_t^{\ast} + \Sigma K_t\bigr)\dd t + \Sigma\,\dd W_t, \\[0.2em]
    \begin{aligned}
    \dd V_t = -\Bigl(&q + (c^{o}_t)^\top u_t^{\ast} + \tfrac12 u_t^{\ast\top}Ru_t^{\ast}- \bigl(\tfrac{\partial V}{\partial Z}\bigr)^{\!\top} K_t\Bigr)\dd t\\
      &+ \bigl(\tfrac{\partial V}{\partial Z}\bigr)^{\!\top}\Sigma\,\dd W_t,
    \end{aligned}
  \end{cases}
\end{align*}
with $V_T = \phi(Z_T)$.
\end{theorem}

Setting $K_t\equiv 0$ recovers the original dynamics; nonzero $K_t$ is used for the importance
sampling of \cite[Section 6]{exarchos2018stochastic}. We use $K_t\equiv0$ in all
experiments reported here.

\section{Deep FBSDE implementation}\label{sec:layer}

This section turns the constrained Hamiltonian minimization into an algorithm: a time discretization and loss for learning the value function and its gradient, the smoothing the second-order conditions require, and the quadratic program solved at each epoch.

\subsection{Discretization}
Following \cite[Section 5.1]{exarchos2018stochastic} we select a time grid
$\{t_0<\dots<t_N=T\}$, take it uniform with $\Delta t = t_{n+1}-t_n = T/N$, and write
$\Delta W_n = W_{t_{n+1}}-W_{t_n}$ for the Brownian increment, so that
$\Delta W_n\sim\mathcal{N}(0,\Delta t\,I_{2M})$. 

The value function is learned in normalized units with a scale $C>0$, and we write $\bar V = V/C$ for the value
of \eqref{eq:cost} divided by it. Every term of the backward equation is divided by $C$, so $\bar V$ at initial time and $\partial \bar V/\partial Z$ is what
the networks represent. The scale is estimated once before training,
as the mean absolute cost-to-go over a pilot batch of rollouts, then held fixed and stored with the
checkpoint. It carries no modeling content, and Section~\ref{sec:layer_qp} shows that it cancels
from the control.

The forward process is discretized by the Euler-Maruyama scheme, initialized by $\tilde V_0 = \bar V_{0}^{\theta_1}$ and
$\tilde Z_0 = Z_0$, where $\bar V_{0}^{\theta_1}(Z)$ approximates the value $\bar V(Z, 0)$ with the neural net having parameter $\theta_1$. The coupled system becomes
\begin{align*}
  \tilde Z_{n+1} =& \tilde Z_n + \bigl(f(\tilde Z_n,Q_{t_n})+ G u^{\ast}_{t_n}(\tilde Z_n)\bigr)\Delta t
    + \Sigma\,\Delta W_n, \\
  \tilde V_{n+1} =& \tilde V_n - \frac{1}{C}\Bigl(q(\tilde Z_n,Q_{t_n})+ (c^{o}_{t_n})^\top u^{\ast}_{t_n}+ \tfrac12 u^{\ast\top}_{t_n} R u^{\ast}_{t_n}\Bigr)\Delta t\notag\\
    &+ \Psi^{\theta_2}(\tilde Z_n,t_n)^{\!\top}\Sigma\,\Delta W_n,
\end{align*}
where $\Psi^{\theta_2}(Z,t)$ approximates the value
gradient $\partial \bar V/\partial Z\in\R^{2M}$ with the neural net having parameter $\theta_2$. 
The loss is the terminal mismatch
\begin{align*}
  \mathcal{L}(\theta_1, \theta_2) = \E_{\theta_1, \theta_2}\bigl[\bigl(\tilde V_N - \phi(\tilde Z_N)/C\bigr)^2\bigr]
    + \lambda_{\mathrm{ridge}}\bigl\|\theta_{2}\bigr\|_2^2 ,
\end{align*}
with $\lambda_{\mathrm{ridge}}>0$, minimized by stochastic gradient descent over sampled initial states and Brownian paths, with the
ridge penalty applied to $\theta_2$.

\subsection{Smoothing of drift}\label{sec:smoothing}

As discussed in Section \ref{sec:cbf}, the unsmoothed drift \eqref{eq:flows} contains a minimum and a positive part, so it is not
differentiable. We replace them with smooth minimum and softplus function
\begin{align*}
  &\mathrm{smin}_{\mu_{\min}}(x,y) = -\tfrac{1}{\mu_{\min}}\log\bigl(e^{-\mu_{\min} x}+e^{-\mu_{\min} y}\bigr),
  \\
  &\mathrm{sp}_{\mu_{+}}(x) = \tfrac{1}{\mu_{+}}\log\bigl(1+e^{\mu_{+} x}\bigr),
\end{align*}
which approximate $\min$ and $(\cdot)^{+}$. We take
$\beta_j,a_j,\beta'_j,a'_j$ to be the derivatives of the composed maps, and the resulting expressions were checked
against automatic differentiation.

\subsection{The safety layer}\label{sec:layer_qp}

At each epoch, and for each sampled trajectory, the order quantity solves
\begin{align}\label{eq:qpfull}
\begin{split}
  &\min_{u,\,s_a,\,s_u\ \ge 0}\,
  \tfrac12 u^\top R_c u + q_{\mathrm{lin}}^\top u
    + \rho\,\bigl(\mathbf{1}^\top s_a + \mathbf{1}^\top s_u\bigr)\\
  &\quad\text{ s.t. }\,\,
   u \ge \underline u, \, u \le \tilde u + s_u,\, W_{\mathrm{agg}} u \le r_{\mathrm{agg}} + s_a .
\end{split}
\end{align}

The objective is the Hamiltonian of \eqref{eq:qp} with the quadratic term scaled by the parameter $\eta>0$, and the constraints are the conditions of
Section~\ref{sec:cbf}, with the following values.
\begin{itemize}
\item $q_{\mathrm{lin}} = G^{\top} \Psi^{\theta_2} + c^{o}_t/C$, the learned marginal value of
      pipeline inventory plus the unit order cost, and $R_c=\eta R/C$, both in the normalized units. 
\item $\underline u = \max\{\underline u_L,\underline u_X,0\}\in\R^M$, the smallest order the
      per-product barriers \eqref{eq:lbL} and \eqref{eq:lbX} admit, and $\tilde u$, the control box
      of \eqref{eq:basic}.
\item $(W_{\mathrm{agg}},r_{\mathrm{agg}})\in\R^{4\times M}\times\R^4$, the four budget rows from
      \eqref{eq:rowL}, \eqref{eq:rowL_2}, \eqref{eq:rowX}, and \eqref{eq:rowX_2}.
\end{itemize}
Every condition is affine in $u$, so \eqref{eq:qpfull} is a quadratic program with $M$ variables, separable across trajectories and warm-started from the previous epoch.

Three devices make it feasible at every solve.
The slacks $s_u\in\R^M$ and $s_a\in\R^4$ relax the box and the budget rows, with penalty rate $\rho = 10\bigl(1+M\|q_{\mathrm{lin}}\|_\infty\bigr)$. The
floor $\beta_{\min}$ bounds $\beta_j$ away from zero, since $\underline u_X$ divides by
$\beta_j$. And $\underline u$ carries no slack, so per-product nonnegativity takes priority over
the box and the budgets whenever they conflict.
Of the back ends CVXPY supports we use CLARABEL, an interior-point method, aligned with a primal-dual interior-point method made by \cite[Algorithm 2]{pereira2021safe}.

\section{Experiments}\label{sec:exp}

We perform three experiments to support the validity of the model. Experiment 1 asks whether the layer holds four aggregate budgets that change on independent piecewise-constant schedules, and how quickly it responds when they do. Experiment 2 isolates what the barrier contributes by removing it in two stages. Experiment 3 measures how performance scales with the assortment size. All three experiments use the instance of Section~\ref{sec:instance} at $M=50$ unless stated, so the three sets of numbers are directly comparable. 

\begin{figure*}[!t]
\centering
\includegraphics[width=0.86\textwidth]{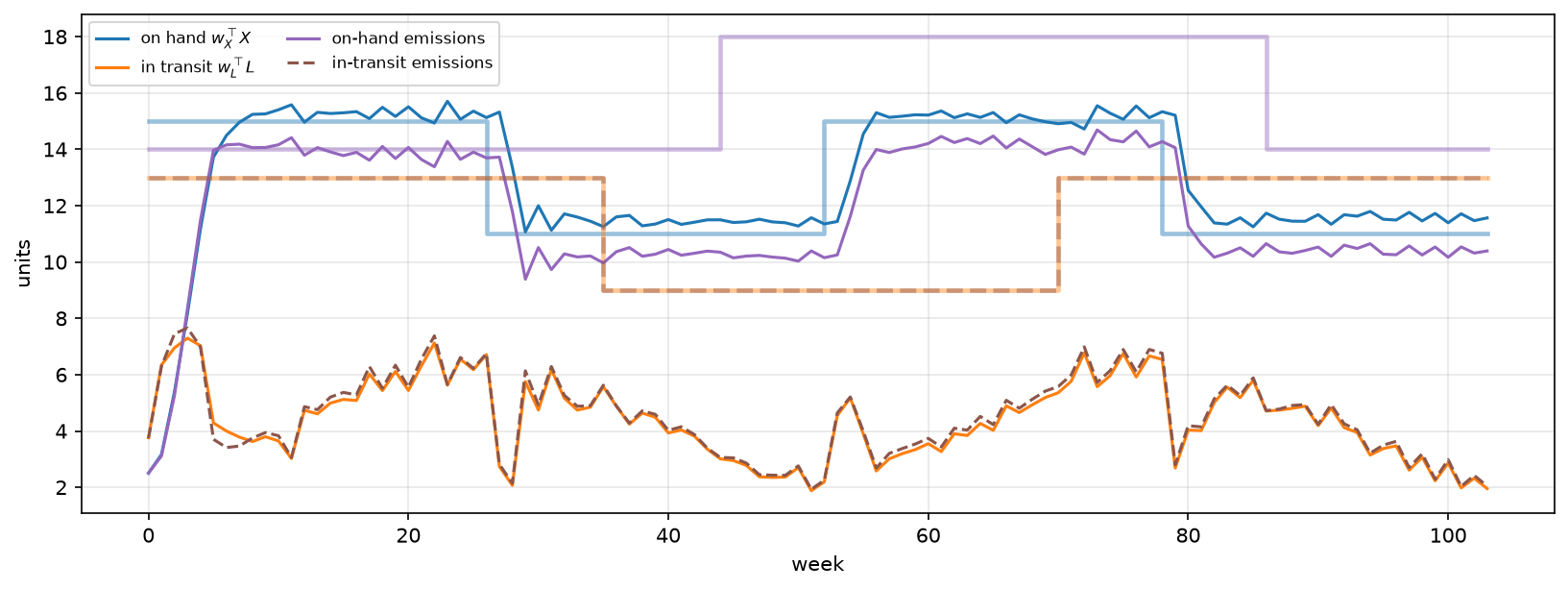}
\caption{The four aggregate budgets of the experimental instance and the aggregates they constrain, over two years of weekly decisions for a $50$-product assortment in Experiment 1. Thick translucent steps are the budgets $H^X$, $H^L$, $H^{E_X}$ and $H^{E_L}$, each moving on its own piecewise-constant schedule; thin curves are the realized aggregates.}
\label{fig:overview}
\end{figure*}

\subsection{Instance}\label{sec:instance}

Time is measured in seasons and every quantity in Table~\ref{tab:instance} is expressed per
season. The horizon is $T=8$, two years of four seasons, discretized into $N=104$ steps so that
one step is one week: $\Delta t=1/13$. The parameter values are generally aligned with \cite{zipkin2008old} for the demand scale and stockout-to-holding
ratio, \cite{maggiar2025structure} for the per-product heterogeneity, and \cite{aviv2001capacitated} for the budget calibration.

Four budgets move on independent schedules, giving nine change events at seven distinct weeks; the emission allowance on in-transit stock tracks $H^L$, so weeks $35$
and $70$ move two rows at once. The sizing condition \eqref{eq:sizing} holds in every one
of the eight resulting intervals.

The order box is $\tilde u(t)=H^L(t)/\Delta t$, which is $9$ or $13$ units per week against a
mean demand of $5$, and is the only control-side limit. It is deliberately generous: it does not
shape normal operation, and its purpose is to keep training away from states in which the value
network was never fitted rather than to represent an operational limit.

\begin{table}[t]
\caption{Instance parameters per season. $U$, Gamma, Exp indicate uniform, Gamma, exponential distribution, respectively.}
\label{tab:instance}
\setlength{\tabcolsep}{2.5pt}
\centering\small
\footnotesize
\begin{tabular}{@{}ll@{}}
\toprule
horizon, grid & $T=8$, $N=104$, $\Delta t=1/13$ \\
demand level, per product & $65\cdot U(0.6,1.4)$ \\
seasonality & $\pm50\%$, one peak per year, trough at $t=0$ \\
demand law & Gamma, coefficient of variation $\sim U(0.2,0.4)$ \\
diffusion & $\sigma_X=\sigma_L=2.7$ \\
holding cost & $c^{h}_j\sim\mathrm{Exp}(13)$ \\
stockout penalty & $c^{p}_j = 9\cdot U(0.6,1.4)$, modulated $\pm 30 \%$\\
order cost & $c^{o}_j = c^{p}_j\cdot(0.2+0.2\,U(0, 1))$  \\
salvage fraction & $\varsigma=0.1$ \\
mean transit time & $\tau_{\mathrm{in}}=1.5/13$ \\
outflow constant & $\tau_{\mathrm{out}}=1/13=\Delta t$ \\
inbound throughput & $g = 2\max_t H^L/\tau_{\mathrm{in}}$\\
capacity weights & $w_X=w_L$ uniform\\
emission weights & $w_{E_X}=w_{E_L}$,$\mathrm{Gamma}(4,4)$ then normalized\\
budgets & $H^X,H^L,H^{E_X},H^{E_L}$ piecewise constant\\
barrier gain & $\gamma=13=1/\tau_{\mathrm{out}}$ \\
regularizer & $R=0.1\,I$, $\eta=8$ \\
control box & $\tilde u(t)=H^L(t)/\Delta t$\\
smoothing, floor & $\mu_{\min}=8/15$, $\mu_{+}=20/15$, $\beta_{\min}=10^{-2}$ \\
\bottomrule
\end{tabular}
\end{table}

\subsection{Metrics}\label{sec:metrics}
In addition to conventional linear cost ($J_{\text{lin}}$), fill rate, and violation, we introduce magnitude measures $M_1, M_3$, largest excess, and utilization for each aggregate constraints. The magnitude measures $M_1, M_3$ are motivated from \cite[Table 1]{eisenach2024neural}, which are normalized because the budget moves.

For example, for the aggregate constraint $w_{X}^{\top} X \leq H^{X}(t), $ we compute $M_1(X)$ and $M_3(X)$ as
\begin{gather}\label{eq:m1m3}
\begin{split}
&M_1(X) \;=\; \E\!\left[\frac{\pos{w_X^\top X - H^X(t)}}{H^X(t)}\right],\\ &M_3(X) \;=\; \Prob\!\left(\frac{\pos{w_X^\top X - H^X(t)}}{H^X(t)} > 0.10\right),
\end{split}
\end{gather}
and we compute $M_1$ and $M_3$ for constraints related to $H^L, H^{E_L}, H^{E_X}$, respectively.

\begin{itemize}
\item $J_{\mathrm{lin}}$: holding $+$ underage $+$ purchase $+$ terminal cost, excluding $R$.
\item fill rate: demand served as a fraction of demand.
\item violation: fraction of (path, week) cells exceeding the row's budget by more than $10^{-6}$.
\item $M_1$, $M_3$: \eqref{eq:m1m3}, mean relative excess; probability of exceeding by more than $10\%$.
\item largest excess: $\max_{t,\text{path}}\pos{w_X^\top X-H^X(t)}/H^X(t)$.
\item utilization: $w_X^\top \cdot X/H^X(t)$; above one places the state outside the safe set.
\item worst per-product position: $\min_{t,j}\bar X^j_t$ where $\bar X^j_t$ is the average on-hand level of product $j$ at week $t$.
\end{itemize}

\subsection{Experiment 1: tracking piecewise-constant budgets}\label{sec:e1}

This experiment measures how closely the dynamics tracks four scheduled budgets, associated to Theorem~\ref{thm:pwc}. Four budgets follow independent
piecewise-constant schedules and the same control law is used throughout. 

From Table \ref{tab:e1rows}, the violation is $0.815$ on $X$ and $0.077$ on $E_X$, and below $10^{-5}$ on $L$ and
$10^{-4}$ on $E_L$. Mean utilization is $1.02$ on $X$, $0.73$ on $E_X$ and $0.38$ on both $L$ and
$E_L$. On the on-hand capacity row, $M_1$ is $0.050$ and $M_3$ is
$0.083$; mean utilization is $1.02$ with a maximum of $1.43$, and the largest excess on any cell is $0.52$. This indicates that in-transit aggregate constraints are well preserved, compared to on-hand aggregate constraints.

From Table \ref{tab:e1seeds}, across the eight seeds the fill rate is $0.562\pm0.031$ on $[0.514,\,0.591]$ and
$J_{\mathrm{lin}}$ is $215{,}009\pm7{,}774$, a spread of $3.6\%$ of the mean. $M_1$ on the on-hand
row is $0.0500\pm0.0021$ on $[0.046,\,0.052]$. The minimum over products and weeks of the path-averaged per-product position is $+0.238\pm0.046$, and $+0.21$ at its lowest across seeds.

\begin{table}[t]
\caption{Experiment 1 by aggregate row, averaged over eight seeds.}
\label{tab:e1rows}
\centering\small
\begin{tabular}{@{}lcccc@{}}
\toprule
& \multicolumn{4}{c}{aggregate row} \\
\cmidrule(l){2-5}
& $X$ (on hand) & $L$ (in transit) & $E_X$ & $E_L$ \\
\midrule
violation & $0.815$ & $0.000$ & $0.077$ & $0.000$ \\
$M_1$              & $0.050$ & $0.000$ & $0.002$ & $0.000$ \\
$M_3$              & $0.083$ & $0.000$ & $0.000$ & $0.000$ \\
largest excess & $0.52$ & $0.00$ & $0.13$ & $0.05$ \\
utilization & $1.02$ & $0.38$ & $0.73$ & $0.38$ \\
\bottomrule
\end{tabular}
\end{table}

\begin{table}[t]
\caption{Experiment 1 across eight seeds at $M=50$, with $512$ evaluation paths.}
\label{tab:e1seeds}
\centering\small
\begin{tabular}{@{}lcc@{}}
\toprule
& mean $\pm$ sd \\
\midrule
fill rate                        & $0.562\pm0.031$  \\
$J_{\mathrm{lin}}$               & $215{,}009\pm7{,}774$ \\
$M_1$, on-hand               & $0.0500\pm0.0021$  \\
worst per-product position & $+0.238\pm0.046$ \\
\bottomrule
\end{tabular}
\end{table}

\begin{figure*}[!t]
\centering
\includegraphics[width=0.86\textwidth]{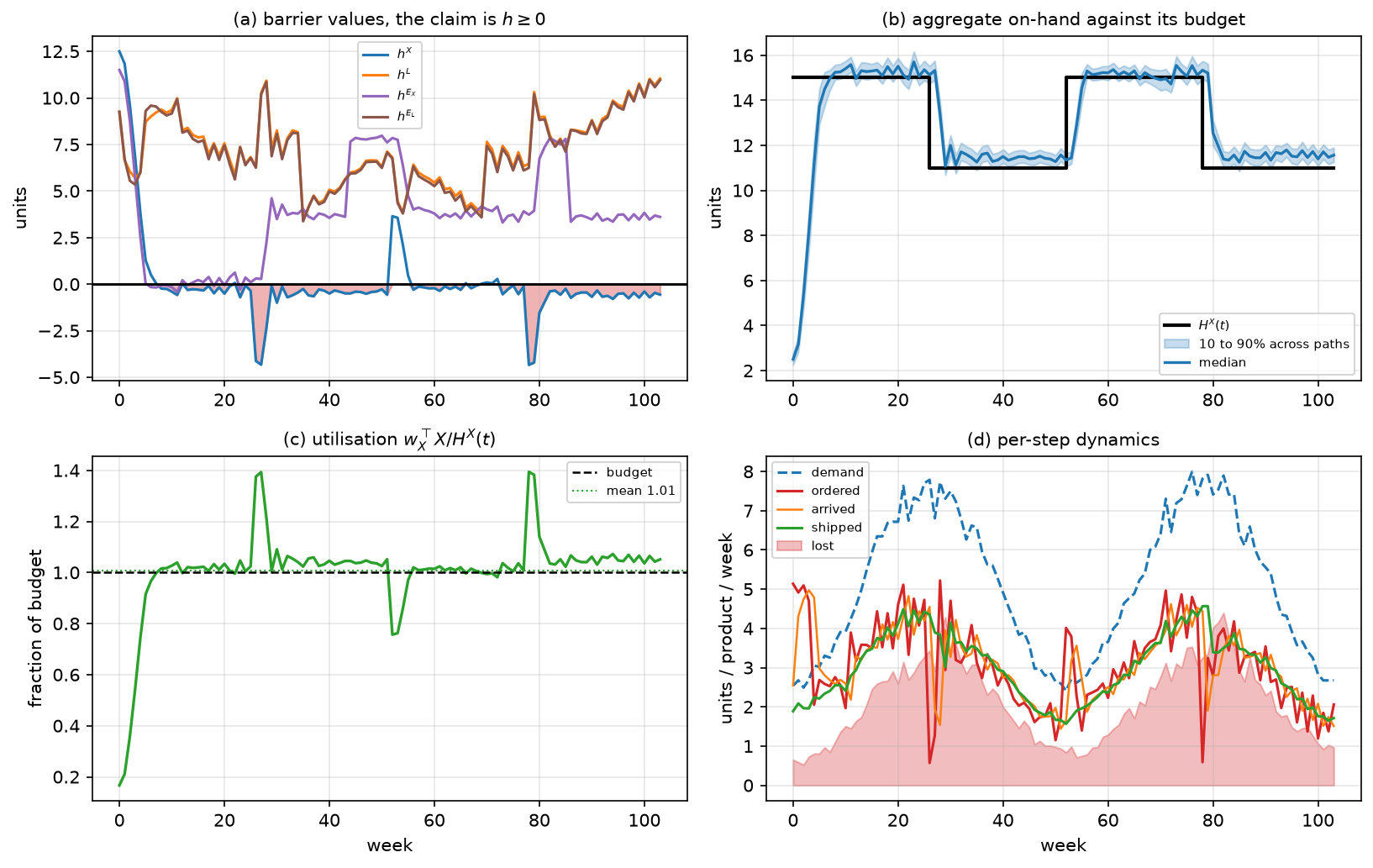}
\caption{Trajectory over the two-year horizon at $M=50$ in Experiment 1. Panel (a) plots the four aggregate barrier
values $h$ against zero. Panel (b) plots the aggregate on-hand row against $H^X(t)$, with a band
from the $10$th to the $90$th percentile across paths. Panel (c) plots utilization $w_X^\top X/H^X(t)$. Panel (d) plots the per-step dynamics: demand, ordered, arrived, shipped and lost.}
\label{fig:behavior}
\end{figure*}

\subsection{Experiment 2: what the safety layer contributes}\label{sec:ablation}

The layer is removed in two stages, giving three arms that differ only in which constraints are
imposed:
\begin{itemize}
\item Full --- the layer of Section~\ref{sec:layer_qp}, per-product barriers and four aggregate budget rows;
\item Relaxed --- the four aggregate budgets raised to $10^4$ so those rows never bind, the per-product barriers retained;
\item Box-only --- as \emph{relaxed}, with $\underline u$ forced to zero, leaving the program $0\le u\le\tilde u$.
\end{itemize}
The control box, instance, network and seeds are identical across arms, and all
three arms are scored against the true schedule. The eight seeds are shared, so comparisons are
paired. The full arm is identical to Experiment 1 by construction and is read off the Experiment 1 runs.
Relaxed and box-only settings differ only in the $2M$ per-product barriers; full and relaxed settings differ only in the four aggregate budget rows.

Between full and relaxed, $M_1(X)$ on the on-hand row is $0.0500$ against $8.69$, and the mean aggregate on-hand is $13.1$ units against $118.7$. Between relaxed and box-only settings, the worst per-product position falls from $+0.159$ to $-5.320$ units. The relaxed arms optimize over a larger feasible set, so their higher fill rate is not directly related to the improvement of the policy.

\begin{figure}[t]
\centering
\includegraphics[width=0.48\textwidth]{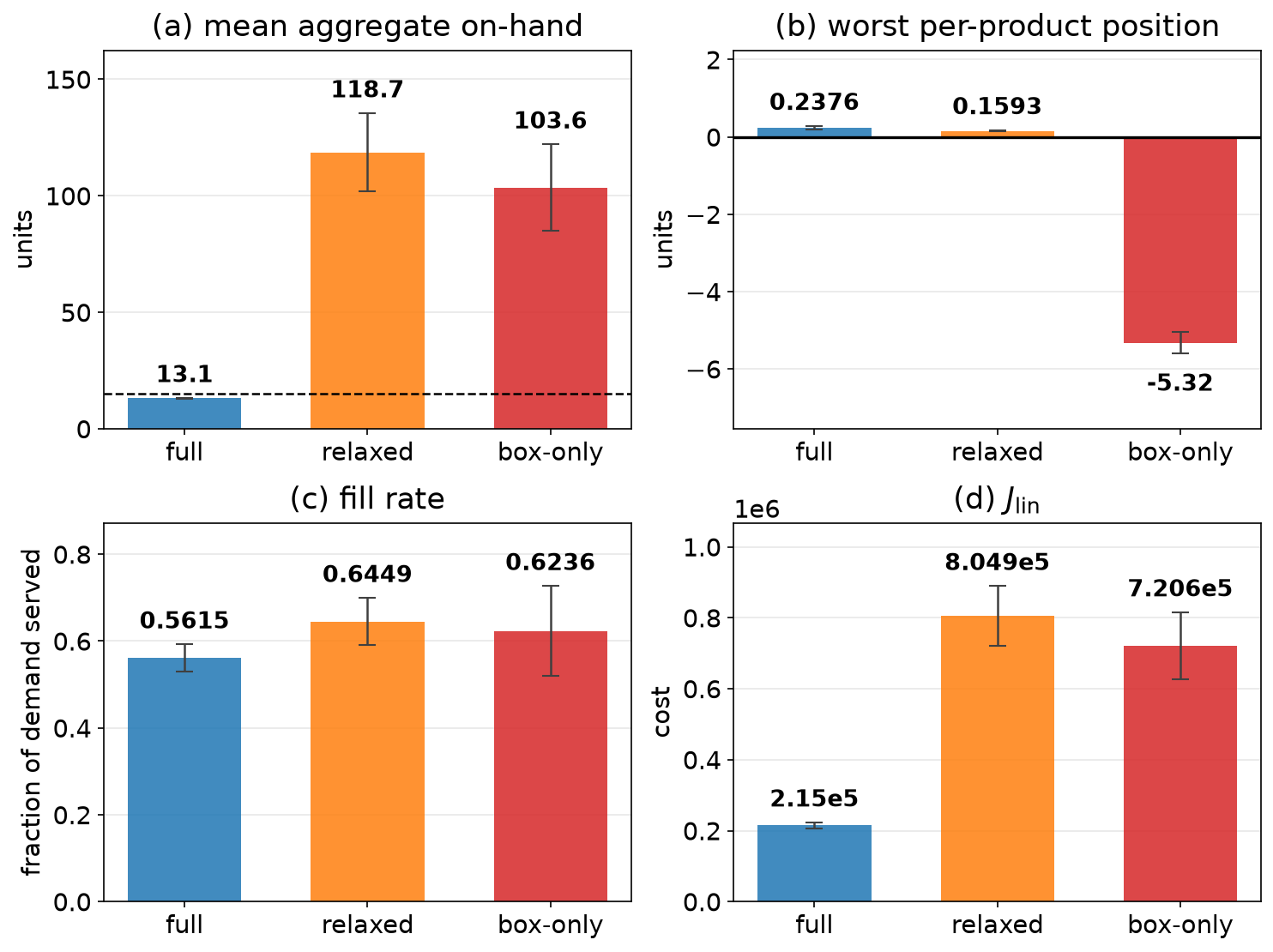}
\caption{The three arms of Experiment 2. Bars are means over eight seeds, error bars the
seed standard deviation. (a) mean aggregate on-hand $w_X^\top X$, with the dashed line at
$\max_t H^X(t)=15$; (b) worst per-product position $\min_{t,j}X^j$, with zero marked; (c) fill rate;
(d) $J_{\mathrm{lin}}$.}
\label{fig:ablation}
\end{figure}

\subsection{Experiment 3: scaling with the assortment}\label{sec:scaling}

We double the assortment from $M=50$ to $M=100$ and scale the network with it across four seeds, increasing the width $w$ as $256$ to $512$. Table~\ref{tab:scaling} reports both sizes on the measures of Tables~\ref{tab:e1rows} and
\ref{tab:e1seeds}. From $M=50$ to $M=100$, $M_1(X)$ is $0.0500$ and $0.0530$, $M_3(X)$
is $0.083$ and $0.072$, and the largest excess is $0.52$ and $0.50$. The violation on the
on-hand row is $0.815$ and $0.875$, on the on-hand emission row $0.078$ and $0.150$, and the two
in-transit rows stay below $10^{-4}$ at both sizes. The fill rate is $0.561$ and $0.601$, and the
worst per-product position $+0.238$ and $+0.314$. Experiment 3 establishes that performance does not degrade when the assortment doubles, when network capacity
is scaled with it.

\begin{table}[t]
\caption{Doubling the assortment on the instance in Experiment 3.}
\label{tab:scaling}
\setlength{\tabcolsep}{2.5pt}
\centering\small
\begin{tabular}{@{}lcc@{}}
\toprule
& $M=50$ & $M=100$ \\
\midrule
fill rate                                  & $0.561\pm0.031$   & $0.601\pm0.021$   \\
$M_1(X)$                       & $0.0500\pm0.0021$ & $0.0530\pm0.0009$ \\
$M_3(X)$                         & $0.083\pm0.005$   & $0.072\pm0.003$   \\
largest excess         & $0.52\pm0.02$     & $0.50\pm0.02$     \\
on-hand violation                & $0.815\pm0.019$   & $0.875\pm0.006$   \\
on-hand emission violation       & $0.078\pm0.054$   & $0.150\pm0.055$   \\
in-transit violation  & $0.000$           & $0.000$           \\
in-transit emission violation  & $0.000$           & $0.000$           \\
worst per-product position & $+0.238\pm0.046$  & $+0.314\pm0.038$  \\
\bottomrule
\end{tabular}
\end{table}

\section{Conclusion}\label{sec:conclusion}

We posed multi-product replenishment under aggregate operating limits as a safe stochastic
optimal control problem, and extended the stochastic high-relative-degree barrier construction
from a single barrier to four coupled aggregate constraints of mixed relative degree, solved as
one quadratic program per sampled trajectory. The construction targets the limits pathwise
rather than in expectation, holds budgets that follow independent piecewise-constant schedules, and retains its performance when the assortment doubles provided
network capacity is scaled with it.

\section*{Acknowledgement}

During the preparation of this work, the authors used Claude Opus 5 (Anthropic) \cite{llmtool} for language editing of the manuscript and for assistance in debugging the implementation code. They reviewed and tested all suggested corrections, verified all reported results, and take full responsibility for the content of this publication.

\bibliographystyle{IEEEtran}
\bibliography{sample}

\end{document}